%
\input amstex
\documentstyle{gen-p}


\def\modo#1{\left|#1\right|}
\def\normo#1{\left\|#1\right\|}
\def\moreproclaim{\par}
\def\conv{\mathop{\hbox{\rm conv}}}
\def\R{{\Bbb R}}

\topmatter
\title 
Random Rearrangements and Operators
\endtitle

\author 
{Stephen Montgomery-Smith and Evgueni Semenov}
\endauthor

\leftheadtext{
Stephen Montgomery-Smith and Evgueni Semenov
}%
\rightheadtext{}%

\address 
Department of Mathematics,
University of Missouri,
Columbia, MO 65211
\endaddress

\email
stephen\@math.missouri.edu
\endemail

\address
Department of Mathematics,
Voronezh State University,
394693 Voronezh, Russia
\endaddress

\email
root\@mathd.vucnit.voronezh.su
\endemail

\issueinfo{00}
{0}
{0}
{1997}

\dedicatory
Dedicated to Professor Selim Krein
on his 80th birthday
\enddedicatory


\endtopmatter

\document

\head{0. Introduction}\endhead
Let $n$ be a positive integer, $x = (x_{ij})_{1\le i,j \le n}$, and
$S_n$ the group of permutations of $\{1,\,2,\dots,n\}$.  Denote
by $x_1^*,\,x_2^*,\dots,x_{n^2}^*$ the rearrangement of $\modo{x_{ij}}$
in decreasing order.  S.~Kwapien and C.~Sch\"utt proved the following
results.

\proclaim{Theorem A {\rm \cite{KS1}}} We have that
$$ {1\over 2n} \sum_{k=1}^n x_k^*
   \le {1\over n!} \sum_{\pi\in S_n} \max_{1\le i \le n} \modo{x_{i,\pi(i)}}
   \le {1\over n} \sum_{k=1}^n x_k^* .$$
\endproclaim

\proclaim{Theorem B {\rm \cite{Sc}}}  If $1\le p \le q < \infty$, then
$$ \eqalignno{
   {1\over 10} & \left(
   \left( {1\over n} \sum_{k=1}^n (x_k^*)^p \right)^{1/p}
   + \left( {1\over n} \sum_{k=n+1}^{n^2} (x_k^*)^q \right)^{1/q}
   \right) \cr
   & \le \left( {1\over n!} \sum_{\pi\in S_n} \left( \sum_{i=1}^n
         \modo{x_{i,\pi(i)}}^q \right)^{p/q} \right)^{1/p} \cr
   & \le
   \left( {1\over n} \sum_{k=1}^n (x_k^*)^p \right)^{1/p}
   + \left( {1\over n} \sum_{k=n+1}^{n^2} (x_k^*)^q \right)^{1/q} .\cr}$$
\endproclaim

There are two ways to generalize these results.  They were presented in
\cite{S1}, \cite{S2} and \cite{M2}.  This article is devoted to the development of these
methods.

Let $1\le q \le \infty$, and $\ell$ be a one-to-one correspondence of
$S_n$ to $\{1,\,2,\dots,n!\}$.  We define the quasi-linear operator
$T_q$ as follows:
$$ T_q x(t) = \left(\sum_{i=1}^n \modo{x_{i,\pi(i)}}^q \right)^{1/q} ,
   \qquad t \in \left[ {\ell(\pi)-1\over n!} , {\ell(\pi)\over n!} \right),$$
with the usual modification for $q = \infty$.  The operator $T_q$ acts from
the set of $n\times n$ matrices into the step functions.  Clearly it depends
upon the choice of $\ell$.  However, if $E$ is a rearrangement invariant
space (see Section~1 for the definitions), then $\normo{T_q x}_E$ does not
depend upon $\ell$.
We define the operator
$$ Ux(t) = \sum_{k=1}^n x^*_k \chi_{((k-1)/n,k/n)}(t) $$
on the set of $n\times n$ matrices $x = (x_{ij})$, where $\chi_e$ is
the characteristic function of $e \subset [0,1]$.

(We would like to mention that there is another modification of this
construction.  If $E$ is an r.i.\ space, we can construct a space
$\tilde E$ on the group $S_n$ equipped with the Haar measure.  In this
case it is not necessary to introduce the function $\ell$.  However, some
additional difficulties appear.)

The inequality
$$ {1\over 12} \left( \normo{U_x}_E +
   \left({1\over n} \sum_{k=n+1}^{n^2} x_k^{*q} \right)^{1/q} \right)
   \le \normo{T_q x}_E
   \eqno (1) $$
was proved \cite{S1} for any matrix $x$, rearrangement invariant space $E$, and
$1 \le q < \infty$.  A more exact estimate is valid for $q = \infty$:
$$ {1\over 2} \normo{Ux}_E \le \normo{T_\infty x}_E \le \normo{Ux}_E .
   \eqno (2) $$
The inverse inequality
$$ \normo{T_q x}_E \le C \left( \normo{Ux}_E +
   \left({1\over n} \sum_{k=n+1}^{n^2} x_k^{*q} \right)^{1/q} \right)
   \eqno (3) $$
was established under the additional assumptions $\alpha_E > 0$, where
$\alpha_E$
is the lower Boyd index of an r.i.\ space $E$, and $C$ does not depend
upon $x$ or $n$.  It is evident that $(3)$ fails for $E = L_\infty$ and
$1 \le q < \infty$.

In this article we generalize these results.  In Chapter~4,
we give a complete criterion for Lorentz spaces for which $(3)$
holds.  In Chapter~5, we consider other r.i.\ spaces for which $(3)$
holds, and give some interpolation results.  All this is based on
Chapter~3, where it is shown that one can reduce $(3)$ to the special
case of diagonal matrices.

In Chapter~7, we consider a completely different generalization, where
in effect $E$ is $L_1$, but $T_q$ is replaced by something analogous to
$T_X$, where $X$ is a symmetric sequence space.  For the analogous result
to $(3)$, one needs the concept of a $D^*$-convex space.  For this reason,
in Chapter~6, we develop the theory of $D$ and $D^*$-convex spaces, building
on earlier work of Kalton \cite{K}.  In Chapter~8 we develop this idea
further, and classify which Lorentz spaces are $D$ or $D^*$-convex.

Some results from this article were announced in \cite{S2} and \cite{M2}.

\head{1. Preliminaries}\endhead
If $x(t)$ is a measurable function on $[0,1]$, we denote by $x^*(t)$
the decreasing rearrangement of $\modo{x(t)}$.  A Banach space $E$
on $[0,1]$ is said to be rearrangement invariant (r.i.) if $y \in E$\
and $x^* \le y^*$ implies that $x\in E$ and $\normo x_E \le \normo y_E$.

The embeddings $L_\infty \subset E \subset L_1$ are true for every r.i.\
space $E$.  In fact, $1 \in E$.  Without loss of generality (except
in Section~7 and parts of Section~6), we
may assume that
$$ \normo 1_E = 1 . \eqno(4)$$
We write $x \prec y$ if
$$ \int_0^\tau x^*(t) \, dt \le \int_0^\tau y^*(t) \, dt $$
for each $\tau \in [0,1]$.  If $E$ is separable or isometric to the
conjugate of some separable r.i.\ space, that $x \prec y$ implies
$\normo x_E \le \normo y_E$.  Denote
$$ E' = \left\{ x : x \in L_1 , \, \int_0^1 xy\, dt < \infty\ \ \forall y \in E
   \right\} ,$$
and equip it with the norm
$$ \normo x_{E'} = \sup_{\normo y_E \le 1} \int_0^1 xy \, dt .$$
Then $E'$ is an r.i.\ space.  The embedding $E \subset E''$ is isometric.
In fact $\normo x_E = \normo x_{E''}$ for all $x \in L_\infty$.
The function $\normo{\chi_{(0,s)}}_E$ is called the fundamental
function of $E$.

Throughout this paper we will assume that all r.i.\ spaces
are either maximal or minimal in the sense of Lindenstrauss
and Tzafriri \cite{LT} (or one may restrict to the case that
the r.i.\ space is separable orisometric to the dual of a separable
r.i.\ space).  Thus for $x\in E$ with $x\ge0$, we have that
$\lim_{t\to \infty} \normo{\min\{x,t\}}_E = \normo x_E$.

Let
$\varphi(t)$ be an increasing function from $[0,1]$ to $[0,1]$,
with $\varphi(0)=0$ and $\varphi(1) = 1$, and continuous on $(0,1]$.
Let $1 \le r < \infty$.
The Lorentz space $\Lambda_r(\varphi)$ consists
of those functions on $[0,1]$ for which the functional
$$ \normo x_{\Lambda_r(\varphi)} =
   \left(\int_0^1 (x^*(t))^r \, d(\varphi(t))^r \right)^{1/r} $$
is finite.

Let us set $\Phi_r$ to be the collection of those $\varphi$ satisfying
the above conditions, and also that $\varphi(t)^r$ is
concave.  Then we see that if $\varphi$ is in $\Phi_r$, then
$\Lambda_r(\varphi)$ satisfies the triangle inequality, and hence is an
r.i.\ space.

Also, for $r>1$,
it is also known that this is equivalent to a norm satisfying the
triangle inequality if and only if there exist $c>0$ and $\epsilon>0$
such that $\varphi(ts) \ge c^{-1} t^{1-\epsilon} \varphi(s)$ for
$0 \le s,t\le 1$ (see \cite{Sa}).

In the earlier sections, we shall primarily be interested in the case when
$r = 1$, and so we will write $\Lambda(\varphi)$ for $\Lambda_1(\varphi)$,
and $\Phi$ for $\Phi_1$.

If $\varphi \in \Phi$ is continuous at $0$, then $(\Lambda(\varphi))^*$
is equal to the Marcinkiewicz space,
$M(\varphi)$, where $M(\varphi)$ is the space of functions on $[0,1]$ for
which the functional
$$ \normo x_{M(\varphi)} = \sup_{0<s\le 1}
   { \int_0^s x^*(t) \, dt \over \varphi(s) } $$
is finite.  This space is also an r.i.\ space.
The spaces $(\Lambda(\varphi))'$ and $M(\varphi)$ coincide for every
$\varphi \in \Phi$.

If $\varphi \in \Phi$, then $\varphi(t)$ and $t/\varphi(t)$ increase on
$[0,1]$.  A function having these properties is said to be quasi-concave.
If $\psi$ is quasi-concave, then there exists $\varphi \in \Phi$ such
that ${1\over2} \varphi \le \psi \le \varphi$.  Indeed, $\varphi$ may
be chosen as the concave majorant of $\psi$.

Let $M(t)$ be a concave even function on $(-\infty,\infty)$, $M(0) = 0$.
Then we define another r.i.\ space,
the Orlicz space $L_M$, to consist of all functions on $[0,1]$ for which
the functional
$$ \normo x_{L_M} = \inf \left\{ \lambda :\ \lambda >0 , \,
   \int_0^1 M \left( {\modo{x(t)} \over \lambda}\right) \, dt \le 1
   \right\} $$
is finite.

We shall use Peetre's $K$-method in this article.  Therefore we
present the main definitions.  Let $(E_0,E_1)$ be a compatible pair
of Banach spaces.  The $K$-functional is defined for each $x \in E_0 + E_1$
and $t>0$ by
$$ K(t,x) = K(t,x,E_0,E_1) = \inf(\normo{x_0}_{E_0} + t\normo{x_1}_{E_1}) ,$$
where the infimum is taken over all representations $x = x_0 + x_1$ with
$x_0 \in E_0$ and $x_1 \in E_1$.  If $0<\theta<1$ and $1 \le p \le \infty$,
then the space $(E_0,E_1)_{\theta,p}$ consists of all $x \in E_0+E_1$
for which the functional
$$ \normo x_{\theta,p} = \left( \int_0^\infty (t^{-\theta} K(t,x))^p
   \, {dt\over t} \right)^{1/p} $$
is finite.  The space $(E_0,E_1)_{\theta,p}$ is an interpolation space
with respect to $E_0,E_1$.  All the above mentioned properties of r.i.\
spaces and the $K$-method can be found in \cite{BS}, \cite{BK}, 
\cite{KPS}, \cite{LT}.

\head{2. Hardy-Littlewood semiordering}\endhead
The semiordering $\prec$ can be applied to vectors.  The definition is
completely analogous.  In this section, we shall establish some
preliminary statements about this semiordering.

If $x,y \in \R^n$, $x,y \ge 0$,
$$ X(t) = \sum_{k=1}^n x_k \chi_{((k-1)/n,k/n)}, \quad
   Y(t) = \sum_{k=1}^n y_k \chi_{((k-1)/n,k/n)}, $$
then the correlations $x \prec y$ and $X \prec Y$ are equivalent.

In general, if $a,b,c \in \R^n$, $a,b,c \ge 0$, then it does not follow
from $a \prec b$ that $a+c \prec b+c$.

Given $x \in \R^n$, $I \subseteq \{1,2,\dots,n\}$, denote
$$ (x \chi_I)_k = \cases
   x_k & \text{ if } k \in I \\
   0 & \text{ if } k \notin I,\  1 \le k \le n. \endcases $$

\proclaim{Lemma 1}  Let $\{1,2,\dots,n\} = I_1 \cup I_2 \cup \dots
\cup I_s$, where $I_1,\,I_2,\dots,I_s$\ are disjoint subsets.
If $x\chi_{I_k} \prec y\chi_{I_k}$ for all $k = 1,\,2,\dots,s$, then
$x \prec y$.
\endproclaim

The proof is obvious.  This statement remains true if we consider
$x\chi_{I_k}$ and $y\chi_{I_k}$ as elements of $\R^{\modo{I_k}}$.

Let $u = (u_{i,j})_{1\le i,j, \le n}$, $u\ge 0$, $u_{11}>0$,
$u_{i,2} = 0$ for $i=1,\,2,\dots,n$.  Put
$$ v_{ij} = \cases
   0 & \text{ if } i=j=1 \\
   u_{11} & \text{ if } i=1,\ j=2 \\
   u_{ij} & \text{ otherwise}.\endcases $$
Given an $n\times n$ matrix $x$, consider the vector
$$ (Tx) = \left(
   \sum_{k=1}^n x_{k,\pi(k)} ,\ \, \pi \in S_n \right) $$
as an element of $\R^{n!}$.

\proclaim{Lemma 2}  $Tu \prec Tv$.
\endproclaim

\demo{Proof}  Denote
$$ \eqalignno{
   S_0 &= \{ \pi\in S_n : \pi(1) \ne 1,\, \pi(2) \ne 1 \} \cr
   S(1,j) &= \{ \pi\in S_n : \pi(1) = 1,\, \pi(2) = j \} \cr
   S(j,1) &= \{ \pi\in S_n : \pi(1) = j,\, \pi(2) = 1 \} ,\cr }$$
where $j = 2,\,3,\dots,n$.  It is evident that
$$ S_n = S_0 \cup \bigcup_{j=2}^n (S(1,j) \cup S(j,1)) \eqno (5)$$
is a disjoint decomposition of $S_n$.  Clearly
$$ (Tu)\chi_{S_0} = (Tv)\chi_{S_0} . \eqno (6)$$
Put
$$ (\tilde T x)(\pi) = x_{1,\pi(1)} + x_{2,\pi(2)} .$$
If $\mu \in S(1,j)$ for some $j\in\{2,\,3,\dots,n\}$, then there exists a
unique $\nu\in S(j,1)$ such that $\mu(k) = \nu(k)$ for $k=3,\,4,\dots,n$.
Therefore,
$$ \eqalign{
   u_{k,\mu(k)} &= v_{k,\nu(k)} \cr
   u_{k,\nu(k)} &= v_{k,\mu(k)} \cr}
   \eqno(7) $$
for every $k=3,\dots,n$.  For such $\mu,\nu$ we have
$$ ((\tilde T u)(\mu),(\tilde T u)(\nu)) \prec
   ((\tilde T v)(\mu),(\tilde T v)(\nu)) .
   \eqno (8) $$
  From $(6)$ and $(7)$ we get
$$ ((T u)(\mu),(T u)(\nu)) \prec
   ((T v)(\mu),(T v)(\nu)) .$$
Given $j \in \{3,\,4,\dots,n\}$, there exist $\mu_i\in S(1,j)$ and
$\nu_i \in S(j,1)$ such that
$$ S(1,j) \cup S(j,1) = \bigcup_{i=1}^{(n-2)!} \{\mu_i,\nu_i\} $$
and
$$ ((T u)(\mu_i),(T u)(\nu_i)) \prec
   ((T v)(\mu_i),(T v)(\nu_i)) $$
for each $i=1,\,2,\dots,(n-2)!$.  By Lemma~1
$$ (Tu)\chi_{S(1,j)} \prec (Tv)\chi_{S(j,1)} \eqno (9) $$
for every $j=2,\,3,\dots,n$.  Taking into account $(5)$, $(6)$, and
$(9)$, and applying Lemma~1 again, we get that
$Tu \prec Tv$.
\qed
\enddemo

Denote by $Q_n$ the set of $(n\times n)$ matrices $(x_{ij})$ such that
$x_{ij} = 0$ or $1$ for every $1 \le i,j \le n$ and
$$ \modo{\{(i,j): x_{ij} = 1\}} =n .$$
The identity matrix is denoted by $I_n$.

\proclaim{Lemma 3}  If $x\in Q_n$, then $Tx \prec TI_n$.
\endproclaim

\demo{Proof}  Clearly $Tx$ and $TI_n$ are equidistributed for every
permutation matrix $x$.  If $x$ is not a permutation matrix, then there
exists a pair of columns or rows such that the first one contains two
or more $1$'s, and the second one contains no $1$'s.  Without loss of
generality we may assume that $x_{11} = x_{21} = 1$, and that
$x_{i,2} = 0$ for each $i=1,\,2,\dots,n$.
Put
$$ y_{ij} = \cases
   0 & \text{ if } i=j=1 \\
   1 & \text{ if } i=1,\ j=2 \\
   x_{ij} & \text{ otherwise}.\endcases $$
By Lemma~2, $Tx \prec Ty$.  If $y$ is a permutation matrix, then the lemma
is proved.  If $y$ is not a permutation matrix, we can use this construction
again.  We obtain a permutation matrix after less than $n$ iterations.
\qed
\enddemo

\proclaim{Lemma 4}  Let $x$ be a $(n\times n)$ matrix such that
$$ \sum_{i,j=1}^n x_{ij} \le n, \quad 0 \le x_{ij} \le 1.
   \eqno (10)$$
Then
$$ Tx \prec TI_n .\eqno(11) $$
\endproclaim

\demo{Proof}  Denote by $P_n$ the set of $(n\times n)$ matrices
satisfying condition $(10)$.  It is evident that the
set of extremal points of $P_n$ coincides with $Q_n$.
Given $j = 1,\,2,\dots,n$, consider the functional
$$ f_j(x) = \max \sum_{\pi\in R} \sum_{k=1}^n x_{k,\pi(k)} ,$$
where the maximum is taken over all subsets $R \subset S_n$ with
$\modo{R} = j$.  The functional $f_j$ is convex.  Therefore
$$ \max_{x\in P_n} f_j(x) = \max_{x \in \mathop{\hbox{\sevenrm ext}} P_n}
   f_j(x) .$$
If $x \in P_n$, we can find $y \in Q_n$ such that
$ f_j(x) \le f_j(y) $.
By Lemma~3, $f_j(y) \le f_j(I_n)$, and consequently,
$f_j(x) \le f_j(I_n)$.   Thus $(11)$ is proved.
\qed
\enddemo

Let $x,f \in L_1$, and $x,f \ge 0$.  It is well known (\cite{KPS}, II.2.2) that
$x \prec f$ implies $x^\alpha \prec f^\alpha$ for $\alpha \ge 1$.  This
is not true in general if $0<\alpha < 1$.  However, under some additional
assumptions on $f$ the implication $x \prec f \Rightarrow x^\alpha \prec
C \, f^\alpha$ is true.  Denote by $O_f$ the set of $x \ge 0$ such that
$x \prec f$.

\proclaim{Lemma 5}  Let $f \ge 0$, $f \in L_1$, $0<\alpha<1$, $C>1$.
Then
$$ x^\alpha \prec C \, f^\alpha \quad \forall x \in O_f \eqno (12) $$
if and only if
$$ \tau^{1-\alpha} \left(\int_0^\tau f(t) \, dt\right)^\alpha
   \le C \, \int_0^\tau f(t)^\alpha \, dt
   \quad \forall \tau \in [0,1] . \eqno (13) $$
\endproclaim

\demo{Proof}  We may assume that $f = f^*$.  Given $\tau\in (0,1]$,
consider the function
$$ x_\tau(t) = \chi_{(0,\tau)}(t){1\over \tau} \int_0^\tau f(s) \, ds .$$
Since $x_\tau \in O_f$, then $(12)$ implies
$$ \int_0^\tau x_\tau(t)^\alpha \, dt
   = \left({1\over \tau} \int_0^\tau f(s) \, ds \right)^\alpha \tau
   \le C \, \int_0^\tau f(t)^\alpha \, dt .$$
It is equivalent to $(13)$.

Let inequality $(13)$ be valid.  By H\"older's inequality
$$ \int_0^\tau x^*(t)^\alpha \, dt
   \le \left( \int_0^\tau x^*(t) \, dt \right)^\alpha \tau^{1-\alpha}
   \le \left(\int_0^\tau f(t) \, dt\right)^\alpha \tau^{1-\alpha}
   \le C \, \int_0^\tau f(t)^\alpha \, dt .$$
\qed
\enddemo

\proclaim{Lemma 6}  The function $f(t) = T_1 I_n(t)$ satisfies $(12)$
for $\alpha>0$, where $C=6$.
\endproclaim

\demo{Proof}  There exists a sequence $1>\tau_1 \ge \tau_2 \ge
\dots \ge \tau_{n+1} = 0$ such that
$$ f^*(t) = \sum_{j=1}^n j \chi_{[\tau_{j+1},\tau_j]}(t) .$$
The function $f$ is closely connected with the classical coincidence
problem.  It is well known (see \cite{W}, 4.9, 10) that
$$ s_j := \tau_j - \tau_{j+1} = {1\over j!} \sum_{k=0}^{n-j}
   {(-1)^k \over k!} .$$
Denote
$$ q_j = \sum_{k=0}^{n-j} {(-1)^k \over k!} .$$
Then $q_{n-1} = 0$, $q_n = 1$ and ${1\over 3} \le q_j \le {1\over 2}$ for
$j=1,\,2,\dots,n-2$.  Hence
$$ \tau_j = \sum_{i=j}^n s_i \le 3 s_j \eqno(14) $$
for $j \ne n-1$.  Since
$$ j s_j = {1\over (j-1)!} \sum_{k=0}^{n-j} {(-1)^k \over k!} ,$$
we get
$$ \sum_{i=j}^n i s_i \le 3j s_j \eqno(15) $$
for $j \ne n-1$.  Using $(14)$ and $(15)$, we have
$$ \tau^{1-\alpha} \left( \sum_{i=j}^n i s_i \right)^\alpha
   \le (3 s_j)^{1-\alpha} (3 j s_j)^\alpha
   = 3j^\alpha s_j < 3 \sum_{i=j}^n i^\alpha s_i $$
for $1 \le j \le n$.  Consequently, the inequality
$$ \tau^{1-\alpha} \left( \int_0^\tau f^*(t) \, dt \right)^\alpha
   \le 3 \int_0^\tau f^*(t)^\alpha \, dt $$
is proved for $\tau = \tau_j$, $1 \le j \le n$.

If $\tau \in [0,\tau_1]$, we can find $1\le j \le n$ and $\lambda \in [0,1]$
such that
$$ \tau = \cases
   \tau_{j+1} + \lambda s_j & \text{ if }\ \tau > \tau_{n-2}\ \text{ or }
\ \tau \le \tau_n\\
   \tau_n + \lambda s_{n-2} & \text{ if }\ 
\tau_n < \tau \le \tau_{n-2}.\endcases$$
It is sufficient to consider only the case $\tau > \tau_{n-2}$.  By
$(14)$ and $(15)$,
$$ \eqalignno{
   \left(\sum_{i=j+1}^n s_i + \lambda s_j\right)^{1-\alpha} &
   \left(\sum_{i=j+1}^k i s_i + \lambda j s_j\right)^\alpha
   \le
   (3s_{j+1} + \lambda s_j)^{1-\alpha}
   (3(j+1)s_{j+1} + \lambda j s_j)^\alpha \cr
   & \le
   3(j+1)^\alpha(s_{j+1} + \lambda s_j) \cr
   & \le
   6 j^\alpha(s_{j+1} + \lambda s_j) \cr
   & \le
   6 \left( \sum_{i=j+1}^n i^\alpha s_i + \lambda j^\alpha s_j \right).\cr}$$
The obtained inequality shows that
$$ \tau^{1-\alpha} \left(\int_0^\tau f^*(t) \, dt \right)^\alpha
   \le 6 \int_0^\tau f^*(t)^\alpha \, dt .$$
\qed
\enddemo

\head{3. Reduction to diagonal matrices}\endhead
Given an integer $n$, denote by $D_n$ the set of diagonal matrices.  It
is evident that if $x \in D_n$, then $x_k^* = 0$ for $n<k\le n^2$.

\proclaim{Theorem 7}  Let $C,q > 1$, and let $E$ be an r.i.\ space.
If
$$ \normo{T_q y}_E \le C \, \normo{U y}_E $$
for any $y \in D_n$, then
$$ \normo{T_q x}_E \le 7C\, \left(
   \normo{U x}_E + \left( {1\over n} \sum_{k=n+1}^{n^2} x_k^{*q} \right)
   ^{1/q} \right) $$
for any $(n\times n)$ matrix $x$.
\endproclaim

\demo{Proof}  Let
$$ \eqalignno{
   \normo{U x}_E &\le 1 & (16) \cr
   {1\over n} \sum_{k=n+1}^{n^2} x_k^{*q} & \le 1 .  & (17) \cr } $$
We can find matrices $y = (y_{ij})$ and $z = (z_{ij})$ such that $x = y+z$,
$\modo{\mathop{\hbox{\rm supp}} y} \le n$, and $U x = U y$.  Denote by
$w$ a diagonal matrix such that $U w = U y$.  By Lemma~2,
$$ T_1 y \prec T_1 w .$$
Hence,
$$ \normo{T_1 y}_E \le \normo{T_1 w}_E \le C \, \normo{U w}_E =
   C \, \normo{U x}_E $$
and
$$ \normo{T_q y}_E \le \normo{T_1 y}_E \le C\,\normo{U x}_E \le C .
   \eqno(18) $$
If $x_{n+1}^* > 1$, then
$$ \min_{0\le t \le 1} U x(t) = x^*_n \ge x_{n+1}^* > 1 .$$
By $(4)$,
$$ \normo{U x}_E > \normo{1}_E = 1.$$
The obtained inequality contradicts $(16)$.  Consequently, for
every $1 \le i,j \le n$
$$ \modo{z_{ij}} \le x_{n+1}^* \le 1 .$$
Denoting $\modo{x_{ij}}^q$ by $v_{ij}$, we get
$$ \eqalignno{
   0 \le v_{ij} & \le 1 , \qquad 1 \le i,j, \le n \cr
   \sum_{i,j = 1}^n v_{ij} & \le n \cr
   \normo{T_q z}_E &= \normo{(T_1 v)^{1/q}}_E .\cr } $$
By Lemma~4,
$$ T_1 v \prec T_1 I_n .$$
Applying Lemmas~5 and~6, we have
$$ (T_1 v)^{1/q} \prec 6 ( T_1 I_n)^{1/q} .$$
We have mentioned in Section~1 that this inequality implies
$$ \normo{(T_1 v)^{1/q}}_E \le 6 \normo{(T_1 I_n)^{1/q}}_E .$$
Since
$$ \normo{(T_1 I_n)^{1/q}}_E = \normo{T_q I_n}_E \le
   C \, \normo{U I_n}_E = C ,$$
then
$$ \normo{T_q z}_E \le 6C .$$
Using the obtained inequality and $(18)$, we get
$$ \normo{T_q x}_E \le \normo{T_q y}_E + \normo{T_q z}_E \le 7C .$$
\qed
\enddemo

Theorem~7 and $(1)$ lead to the following.

\proclaim{Corollary 8}  Let $E$ be an r.i.\ space, and $1 \le q < \infty$.
The equivalence
$$ \normo{T_q x}_E \approx
   \normo{U x}_E + \left( {1\over n} \sum_{k=n+1}^{n^2} x_k^{*q}
   \right)^{1/q} , \eqno(19) $$
where the equivalence constants depend neither upon the matrix
$x$ nor on $n$, takes place if and only if the estimate
$$ \normo{T_q y}_E \le C \, \normo{U y}_E $$
is valid for every diagonal matrix $y$.
\endproclaim

Denote by $\normo{T_q}_E$ the least $C$ in the last inequality, and by
$F_q$ the set of r.i.\ spaces satisfying condition $(19)$.  Given an r.i.\
space $E$, denote by $\omega(E)$ the set of $q \in [1,\infty]$ such that
equivalence $(19)$ takes place, and put
$$ \tau(E) = \inf \omega(E) .$$
The monotonicity of the function $q \mapsto \normo{T_q x}_E$, $(21)$, and
Corollary~8 imply that $w(E) = [\tau(E),\infty]$ or $\omega(E) = (\tau(E),
\infty]$.  Some examples show that both of these possibilities may be
realized.

\head{4. Lorentz spaces}\endhead
Given $0<j \le k < n$, denote
$$ I_{n,k} = \mathop{\hbox{\rm diag}}(
\underbrace{1,1,\dots,1}_k,\underbrace{0,0,\dots,0}_{n-k}) $$
and put
$$ \mu_{n,k,j}=\mathop{\hbox{\rm mes}} \{ t\in[0,1]:T_1 I_{n,k}(t) = j \}.$$

\proclaim{Lemma 9}  If $0<s<1$, then
$$ {s^j\over e j!} \le \sup \mu_{n,k,j} \le {s^j\over j!} , \eqno(20)$$
where the supremum is taken over $(n,k)$ such that $k \le n s$.
\endproclaim

\demo{Proof}  We have
$$ \eqalign{
   \mu_{n,k,j}
   &\le {C_k^j (n-j)! \over n!}
   = {k!(n-j)!\over j!(k-j)!n!} \cr
   &= {k(k-1)\dots(k-j+1) \over j! n(n-1)\dots(n-j+1)}
   \le {1\over j!} \left({k\over n}\right)^j
   \le {s^j \over j!}. \cr} \eqno(21) $$
Following (\cite{W}, 4.9.B), denote $B_{n,k,j} = n! \mu_{n,k,j}$.  It is known
that
$$ B_{n,k,j} = {j+1 \over k+1} B_{n+1,k+1,j+1} .$$
Therefore,
$$ \mu_{n,k,j} = {1\over n!} B_{n,k,j}
   = {k \over n! j} B_{n-1,k-1,j-1} = {k \over n j} \mu_{n-1,k-1,j-1} .$$
Hence
$$ \eqalignno{
   \mu_{n,k,j}
   &= {k(k-1)\dots(k-j+1)\over j! n(n-1)\dots(n-j+1)} \mu_{n-j,k-j,0} \cr
   &= {k(k-1)\dots(k-j+1)\over j! n(n-1)\dots(n-j+1)}
   \left(1-{1\over n-j}\right)^{k-j} .\cr}$$
The last part tends to $s^j e^{-s}/j!$ if $k=[ns]$ and $n$ tends to infinity.
This proves the left part of $(20)$.
\qed
\enddemo

The above proved statement allows us to solve completely the problem
on the validity of equivalence (19) in the class of Lorentz spaces.
Recall that we denote by $\Phi$ the set of increasing concave functions
on $[0,1]$ with $\varphi(0)=0$ and $\varphi(1) = 1$.

\proclaim{Theorem 10}
Let $\varphi \in \Phi$, $1\leq q<\infty$. The equivalence
$$
\|T_q x\|_{\Lambda (\varphi)} \approx \|U x\|_{\Lambda(\varphi)} +
\left({ 1\over  n} \sum_{k=n+1}^{n^2} {x_{k}}^q\right)^{1/q}
\eqno (22) $$
takes place if and only if
$$
\Gamma_{\varphi,q}:= \sup_{0< t \leq 1} {1\over \varphi(t)} \sum_
{j=1}^{\infty} j^{(1/ q)-1} \varphi({t^j}/{j!}) < \infty.
\eqno (23) $$
Moreover,
$$
\|T_q \|_{\Lambda(\varphi)} \leq \Gamma_{\varphi, q}
\leq qe\|T_q\|_{\Lambda(\varphi)}.
$$
\endproclaim

\demo{Proof}
We use the notation of Lemma~9. If (22) is fulfilled,
then there exists a constant $C>0$ such that
$$
\| T_q I_{n,k}\|_{\Lambda(\varphi)} \leq C\|U I_{n,k}\|_{\Lambda(\varphi)}
= C\varphi(k/n).
$$
Since
$$\align
\|T_q I_{n,k}\|_{\Lambda(\varphi)} &= \left\|\sum_{j=1}^{k} (j^{1/q}
-(j-1)^{1/q})\chi_{(0,\mu_{n,k,j})} \right\|_{\Lambda(\varphi)} \\
&=
\sum_{j=1}^{k} (j^{1/q} - (j-1)^{1/q})\varphi(\mu_{n,k,j}),
\tag24
\endalign
$$
then it follows by Lemma~9 that for each $t\in(0,1)$ and for
each integer $m$ that we have
$$
\sum_{j=1}^{m} (j^{1/ q } - (j-1)^{1/ q }) \varphi
({t^j} /{ej!}) \leq C\varphi(t).
$$
Hence,
$$
\sum_{j=1}^{\infty} j^{(1/q)-1} \varphi({t^j}/{j!}) \leq
Cqe\varphi(t).
$$
This proves the first part of the theorem.

Now suppose that, for every $t\in (0,1]$,
$$
\sum_{j=1}^{\infty} t^{(1/ q) -1} \varphi({t^j}/{j!})
\leq C\varphi(t).
$$
Applying the obvious inequality
$$
j^{1/ q }-(j-1)^{1/ q} \leq j^{(1/ q) -1}
$$
and (24), we get
$$ \eqalignno{
   \|T_qI_{n,k}\|_{\Lambda(\varphi)}
   &= \sum_{j=1}^{k}
   (j^{1/ q} - (j-1)^{1/ q} \varphi (\mu_{n,k,j}) \cr
   &\leq
   \sum_{j=1}^{k} j^{ (1/ q)-1} \varphi (1/ {j!}
   (k/ n )^j) \cr
   &\leq
   C\varphi(k/ n) \cr
   &=
   C\|U I_{n,k}\|_{\Lambda(\varphi)}. & (25) \cr } $$
All Lorentz spaces $\Lambda(\varphi)$ have the following property
(\cite{KPS}, II.5.2). If a convex functional is uniformly bounded
on the set of characteristic functions, then it is
uniformly bounded on the set of step functions. We apply
this property to matrices. Then (25) implies that
$$
\|T_q y\|_{\Lambda(\varphi)} \leq C\|U y\|_{\Lambda(\varphi)}
$$
for each $y\in D$. By Corollary 8, (22) is valid.
\qed
\enddemo

In other words, $\Lambda(\varphi) \in F_q$ if and only if
$\Gamma_{\varphi,q} <
\infty$. We mention that
$$
\|U x\|_{\Lambda(\varphi)} =\sum_{k=1}^{n} x_{k} (\varphi( k/ n)
-\varphi((k-1)/n)).
$$

Let us study condition (23) in detail.

\proclaim{Lemma 11}
Let $0 < \alpha \leq 1$, $a,q\geq 1$, and $\varphi \in \Phi$ with
$\varphi(t) \leq at^\alpha$ for every $t\in [0,1]$. Then $\Lambda
(\varphi) \in F_q$ and $\Gamma_{\varphi,q} \leq {5a}/\alpha$.
\endproclaim

\demo{Proof}
Since $\Gamma_{\varphi, q} \leq \Gamma _{\varphi,1}$, we
shall estimate only $\Gamma_{\varphi, 1}$.
Given $s\in (0,a^{- 1/ \alpha})$,
we construct the function
$$
\varphi_s(t)=\cases
at^\alpha & \text{ if }\ 0\leq t\leq s,\\
as^\alpha & \text{ if }\ s\leq t \leq as^\alpha,\\
t         & \text{ if }\ as^\alpha \leq t\leq 1.
\endcases
$$
The set of the quasi-concave functions $\varphi_s$ possesses the following
property. If $\varphi\in \Phi$ and $t_1\in (0,1)$, then we can find
$s\in (0,a^{-1/ \alpha})$ such that $\varphi_s (t) \leq \varphi(t)$
for $t\in [t_1,1]$, and $\varphi_s(t) \geq \varphi(t)$ for $t\in [0,t_1]$.
Therefore it is sufficient to obtain the needed estimate only for
the function $\varphi_s$ and $t= as^\alpha$. Put
$$
N= \max\left\{j:\ {(as^\alpha)^j\over j!} \geq s\right\} =
\max\{j:\ (j!)^{1/ j} \leq as^{\alpha-1/ j}\}.
$$
Then
$$
\sum_{j=1}^{\infty} \varphi_s \left({(as^\alpha)^j\over j!}\right)
= Nas^\alpha +
a\sum_{j=N+1}^{\infty} \left({(as^\alpha)^j\over j!}\right)^\alpha \leq
as^\alpha \left(N + a\sum_{j=1}^{\infty} ({j!})^{-\alpha} \right).
$$
Since $(j /3)^j \leq j!$, it follows that
$$
N\leq
\max\{j:\ j\leq 3as^{\alpha-1/ j}\}
\leq
\max (3a, 1/ \alpha).
$$
Hence,
$$
\Gamma _{\varphi,1}\leq \max(3a, 1/ \alpha) + a
\sum_{j=1}^{\infty} (j!)^{-\alpha} \leq
{3a}/\alpha + a \sum_{j=1}^{\infty} 2^{-j\alpha} \leq
{5a} / \alpha.
$$
\qed
\enddemo

The assumption $\varphi(1) = 1$ is essential in Lemma~11. Indeed,
if we let $\psi_\epsilon (t) = \min (\sqrt t, \epsilon)$, then $\Gamma_
{\psi_\epsilon, 1}$ tends to $\infty$ when $\epsilon$ tends to $0$.

\head{5. Interpolation Spaces}\endhead
Theorem~10 may be extended on a wider class of r.i.\ spaces.
Given numbers $\alpha\in (0,1]$, $a\geq 1$, denote by $\Phi(a, \alpha)$
the set of functions $\varphi\in \Phi$ such that $\varphi(t) \leq at^\alpha$
for each $t\in [0,1]$. Let $E$ be an r.i. space. Given $y\in E^\prime$,
$\|y\|_{E^\prime} =1$, we put
$$ \psi_y(t) =
   {\displaystyle \int_{0}^{t} y^*(s) ds + t
    \over \|y\|_{L_1} +1} $$
and
$$
\|x\|_1 = \sup_{\|y\|_{E^\prime}= 1}\|x\|_{\Lambda(\psi_y)}.
$$
We mention that $\psi_y\in \Phi$.

\proclaim{Lemma 12}
\roster\widestnumber\item{ii)}
\item"{i)}"
The norms $\|\cdot\|_{E^{\prime\prime}}$ and $\|\cdot\|_1$
are equivalent and
$$
{1\over 2} \|x\|_E
\leq \|x\|_1 \leq
2 \|x\|_E \qquad \hbox{for all $x\in L_\infty$}.
\eqno(26) $$
\item"{ii)}"
If $L_p\subset E$ for some $p <\infty$ and
$$
\|x\|_E \leq a\|x\|_{L_p} \qquad \hbox{for all $x\in L_p$},
\eqno(27) $$
then $\psi_y\in\Phi(a + 1, 1/ p)$ for every
$y\in E^\prime$ with $\|y\|_{E^\prime} =1$.
\endroster
\endproclaim

\demo{Proof}
i) By the Hardy-Littlewood theorem on rearrangements
(see \cite{KPS}, II.2.2.17) 
it follows that for $x\in E^{\prime\prime}$ we have
$$
\|x\|_{E^{\prime\prime}}=
\sup_{\|y\|_{E^\prime}=1} \int_{0}^{1}x(t)y(t) dt =
\sup_{\|y\|_{E^\prime}=1} \int_{0}^{1}x^*(t)y^*(t) dt =
\sup_{\|y\|_{E^\prime}=1} \|x\|_{\Lambda(\varphi_y)},
$$
where $\varphi_y (t) = \int_{0}^{t}y^*(s)ds$. Since
$$
\|y\|_{L_1} + 1 \leq \|y\|_{E^\prime} + 1 = 2,
$$
then $\psi_y \geq {1\over2} \varphi_y$ and
$$
\|x\|_1 \geq {1\over 2} \|x\|_{E^{\prime\prime}} \qquad \hbox{for all $x\in
E^{\prime\prime}$}.
\eqno(28) $$
On the other hand,
$$
\|x\|_1\leq \|x\|_{E^{\prime\prime}} + \|x\| _{L_1} \leq 2\|x\|_
{E^{\prime\prime}} \qquad \hbox{for all $x\in E^{\prime\prime}$}.
\eqno(29) $$
The norms $\|\cdot\|_E$ and $\|x\|_{E^{\prime\prime}}$ coincide
on $L_\infty$ (\cite{BS}, 1.2.7).
Therefore (28) and (29) imply (26).

ii) If ${1\over p} + {1\over {p^\prime}} = 1$, then
$E^\prime \subset L_{p^\prime}$ and
$$
\|y\|_{L_{p^\prime}} \leq a \|y\|_{E^\prime} = a.
$$
By H\"older's inequality,
$$
\psi_y (t) \leq \int_{0}^{t} y^* (s) ds + t \leq \|y\|_{E^\prime}
t^{1/ p} + t \leq a t^{ 1/ p} + t \leq (a+1) t^{1/ p}.
$$
\qed
\enddemo

\proclaim{Theorem 13} Let $E$ be an r.i. space, $E\supset L_p$ for
some $p < \infty$, and $1\leq q < \infty$. Then (19) is fulfilled, that
is,
$E\in F_q$.
\endproclaim

\demo{Proof}
There is a constant $a > 1$ such that (27) is valid.
By Lemma~12~(ii), $\psi_y \in \Phi(a+1,1/ p)$ for every
$y\in E^\prime$ with $\|y\|_{E^\prime} = 1$. By Lemma~11, we
know that
$$
\Gamma_{\psi_y,q} \leq 5(a +1) p.
$$
Applying the second part of Theorem~10, we get that
$$
\|T_q x\|_{\Lambda (\psi_y)}\leq 5(a+1) p \|U x\| _{\Lambda(\psi_y)}
$$
for each $x\in D$, and $y\in E^\prime$ with $\|y\|_{E^\prime} = 1$. Hence,
$$
\|T_q x\|_1 \leq 5(a+1) p\|U x\|_1.
$$
This and (26) imply that
$$
\|T_q x\|_E \leq 20 (a+1) p\|U x\|_E.
$$
By Corollary~8, it follows that (19) is fulfilled.
\qed
\enddemo

We mention that the conditions
\roster\widestnumber\item{2)}
\item"{1)}" $E\supset L_p$ for some $p <\infty$;
\item"{2)}" $t^{-\alpha} \in E$ for some $\alpha >0$
\endroster
are equivalent.

\bigskip

Let $E_1$, $E_2$ be r.i.\ spaces, $1\leq q < \infty$,
$E_1\supset E_2$ and $E_2 \in F_q$. Does it follow that $E_1 \in F_q$?
Theorem~13 shows that the answer to this question is positive
if $E_2 \supset L_p$ for some $p < \infty$. In general, the
answer is negative. We now show this.

\proclaim{Theorem 14} Let $\varphi\in\Phi$ and $1\leq q < \infty$. The
following conditions are equivalent:
\roster\widestnumber\item{iii)}
\item"{i)}" if $\psi \in \Phi$ and $\psi\leq \varphi$, then
$\Gamma_{\psi,q} <\infty$;
\item"{ii)}" there exists a constant $C>0$ such that if $\psi\in\Phi$ and
$\psi\leq\varphi$, then $\Gamma_{\psi,q} \leq C$.
\item"{iii)}" there exist numbers $\alpha\in (0,1]$ and $a\geq 1$ such
that $\varphi (t) \leq a t^\alpha$ for every $t\in [0,1]$.
\endroster
\endproclaim

\demo{Proof}
The implication (iii) $\Rightarrow$ (ii) was proved in
Lemma~11. The implication (ii) $\Rightarrow$ (i) is trivial. Therefore,
we need prove only the implication (i) $\Rightarrow$ (iii). Let
$$
\varphi \notin \bigcup_{0<\alpha \leq 1 \leq a}\Phi (a,\alpha),
$$
so that
$$
\sup_{0< t \leq 1} \varphi(t) t^{-1/ n } = \infty
\qquad n=1,2,\dots
\eqno (30) $$
Using (30) we can find a sequence $t_n\downarrow 0$ such that
$t_1 =1 $ and
$$
\varphi(t_n) \geq {n\varphi(t_{n-1}) \over t_{n-1}} t_n ^{1/n }
$$
for every $n=2,3\dots\ $. Then
$$
\left(
{t_{n-1} \over \varphi(t_{n-1})}\varphi(t_n) \right)^n
\geq n^n t_n > n! t_n.
\eqno (31) $$
If we put
$$
s_n = {t_{n-1} \over \varphi(t_{n-1})} \varphi(t_n),
$$
then $t_n < s_n < t_{n-1}$ for every $n=2,3,\dots\ $. Define
$$
\psi(t) =\cases
\varphi(t_n)                 & \text{ if }\ t_n \leq t \leq s_n,\\
t \varphi(t_{n-1})/{t_{n-1}} & \text{ if }\
 s_n \leq t\leq t_{n-1},\ n=2,3, \dots,\\
0                            & \text{ if }\ t =0.\endcases $$
The function $\psi(t)$ is quasi-concave on $[0,1]$, and $\psi \leq \varphi$.
By (31) we see that
$$
\psi({s_n^n}/{n!}) = \psi(t_n) = \varphi(t_n) = \psi (s_n).
$$
For every integer $n$ we get
$$
\Gamma_{\psi,q}
\geq
{1\over \psi(s_n)}
\sum_{k=1}^{n}
k^{(1/q)-1} \psi ({s_n^k}/{k!}) =
\sum_{k=1}^{n}
k^{(1/q)-1} >
\sum_{k=1}^{n} k^{-1}.
$$
Hence, $\Gamma_{\psi, q}= \infty$. Denote by $\nu(t)$ the concave majorant
of $\psi (t)$. Then $\nu\in \Phi$, $\nu\leq \varphi$ and $\Gamma_{\nu, q} =
\infty$.
\qed
\enddemo

Theorem~13 is practically an interpolation theorem. It shows that if $E$
is an r.i.\ space, $E= E^{\prime\prime}$ and $E\supset L_p$ for some
$p<\infty$, then $E$ is an interpolation space with respect to the set
$\{\Lambda(\varphi),\ \Gamma_{\varphi,1} < \infty\}$. One can prove that the
assumption $E= E^{\prime\prime}$ may be replaced with the separability
of $E$. Using the $K$-method, we can obtain another sufficient condition
for $E\in F_q$.

\proclaim{Theorem 15} Let $\varphi_0, \varphi_1\in \Phi$ and
$0<\gamma,\theta <1$, and suppose that
the function $\varphi_0^\gamma (t) / \varphi_1 (t)$ increases on $(0,1]$.
Then
$$
(\Lambda (\varphi_0), \Lambda (\varphi_1))_{\theta,\infty}
\approx M(\tilde
\varphi_\theta),
\eqno (32) $$
where
$$
\tilde \varphi_\theta (t) = {t  \over \varphi_\theta (t) } =
{t\over \varphi_0^{1-\theta} (t) \varphi_1^\theta (t)}.
$$
\endproclaim

Theorem~15 is a special case of a more general result which is
contained in \cite{S3}. Using some results on the stability of the
interpolation functions \cite{A}, one can obtain a similar statement.

\proclaim{Lemma 16}
Let $E_0$, $E_1 \in F_q$ be r.i.\ spaces, where $q\in [1, \infty)$,
and suppose that
$E_2$ is an interpolation space with respect to
$E_0,E_1$. Then $E_2 \in F_q$.
\endproclaim

\demo{Proof}
Given an integer $n$, we consider the operator
$$
B_n x = \mathop{\hbox{\rm diag}}
\left( \textstyle n \int_{(k-1)/n}^{k/ n} x(s) \, ds,\
1\leq k\leq n\right).
$$
The operator $B_n$ acts from $L_1$ into the set of diagonal
matrices. The operator $U B_n$ is an averaging operator, and
$U B_n x$ is the conditional expectation of $x$ with
respect to the set of intervals $\{({k-1\over n}, {k\over n}),\
1\leq k \leq n\}$. By Theorem 2.a.4 \cite{LT}, it follows that
$$
\|U B_n\|_E = 1.
$$
Corollary~8 shows that $E\in F_q$ if and only if
$$
\sup_{n} \|T_q B_n\|_E < \infty.
$$
This proves the Lemma.
\qed
\enddemo

\proclaim{Lemma 17} Suppose that $1\leq q < \infty$,
that $\varphi\in \Phi$ satisfies condition~(23),
that $1 < \mu < \lambda$, and that $\varphi^\lambda \in\Phi$.
Then $M(\tilde \varphi^\mu) \in
F_q$.
\endproclaim

\demo{Proof}
Applying (23) and Jensen's inequality, we get that
$$
\sum_{j=1}^{\infty}
j^{(1/q) - 1} \varphi^\lambda ({t^j}/{j!})
\leq
\left(\sum_{j=1}^{\infty} j^{(1/ q) -1} \varphi({t^j}/{j!})\right)^\lambda
\leq C^\lambda \varphi^\lambda (t)
$$
for every $t\in (0,1]$. It means that $\varphi^\lambda$ satisfies condition
(23) with the constant $C^\lambda$. By Theorem~15,
$M(\tilde \varphi^\mu)$ is
an interpolation space with respect to $\Lambda(\varphi)$ and
$\Lambda (\varphi^\lambda)$.
The required statements now follow from Lemma~16.
\qed
\enddemo

Consider the following example. Given $p> 0$, we put
$$
\varphi_p (t) = \left({\log(1 + 1/t)\over \log 2}\right)^{-1/p}.
\eqno (33) $$
If $p\geq 1$, then $\varphi_p\in \Phi$. If $p\in (0,1)$, then
$\varphi_p(t)$ is concave in a sufficiently small neighborhood of
the origin. Consequently, $\varphi_p$ is concave up to equivalence.

\proclaim{Lemma 18} Let $p>0$ and $1\leq q<\infty$.
Then $\Gamma_{\varphi_p,q}
<\infty$ if $p<q$. If $\varphi\in \Phi$ and $\varphi\geq C\varphi_q$ for
some $C>0$, then $\Gamma_{\varphi, q} = \infty$.
\endproclaim

\demo{Proof}
Let $p<q$ and $0<t\leq 1$. We have that
$$ \eqalignno{
   {1\over \varphi_p(t)}
   \sum_{j=1}^{\infty}
   j^{(1/ q) -1} \varphi_p ({t^j}/{j!})
   &\leq
   1 + \sum_{j=2}^{\infty}
   j^{(1/ q) -1}
   \left( {\log ( 1 + 1/t) \over \log (1 + {j!}/{t^j})}\right)^{1/p} \cr
   &\leq
   1 + \sum_{j=2}^{\infty}
   j^{(1/ q) -1} \left({1 \over j-1}\right) ^{1/p}
   < \infty. \cr }
$$
Therefore $\Gamma_{\varphi_p, q} < \infty$. If $\varphi \geq C\varphi_q$ and
$0<t\leq 1$, then
$$ \eqalignno{
   \sum_{j=2}^{\infty} j^{(1/ q) -1} \varphi ({t^j}/{j!})
   &\geq
   C(\log 2)^{-1/q } \sum_{j=2}^{\infty} j^{(1/ q) -1}
   \log^{-1/q}(1 + {j!}/{t^j}) \cr
   &\geq
   C\sum_{j=2}^{\infty} j^{(1/ q)-1} (j\,\log(j/t))^{-1/q} \cr
   &=
   C\sum_{j=2}^{\infty} j^{-1} (\log(j/t))^{-1/q} = \infty. \cr }$$
\qed
\enddemo

Let us consider the Orlicz space $\exp L_p$. It is generated by the
function
$$
M_p (u) = e^{|u|^p}-1.
$$
If $p\geq 1$, then $M_p (u)$ is convex and the fundamental function
of $\exp L_p$ is equal to $(\log(1 + 1/t))^{-1/p}$. If
$0 < p < 1$, then $M_p (u)$ is convex for sufficiently large $u$.

\proclaim{Theorem 19} Let $1\leq q < \infty$. If $p<q$, then $\exp L_p
\in F_q$. If $p> q$, then $\exp L_p \notin F_q$.
\endproclaim

\demo{Proof}
Let $p\geq 1$. By \cite{Lo}, the spaces $\exp L_p$ and
$M(\tilde\varphi_p)$ coincide, where $\varphi_p$ is defined by (33).
Applying Lemma~18, we get that $\exp L_p\in F_q$ if $p<q$. The first
part of the theorem is proved.

If $E$ is an r.i.\ space and $E\in F_q$, then
$$
\sup_{n} \|T_qI_n\|_E <\infty.
$$
In fact, let $E$ be an Orlicz space $L_M$. Lemma~9 shows that
$$
\sum_{j=1}^{\infty} { e^{\epsilon^p j ^{p/ q}} - 1 \over j!}
< \infty
$$
for some $\epsilon > 0$.
This series diverges for any $p> q$ and $\epsilon > 0$. Hence,
$p\leq q$.

So, the theorem has been proved for $p\geq 1$. If $0 < p < 1$,
we can change $M_p(u)$ for a convex equivalent function.
Therefore the theorem is valid for every $p>0$.
\qed
\enddemo

\head{6. $D$ and $D^*$-convex Spaces}\endhead
The notion of $D$-convexity was introduced by Kalton \cite{K} (Section~5).
Indeed, much of the proof of this section is inspired by his proof of
Lemma~5.5.

Given a function $x$ on $[0,1]$, we will define its distribution
function $d_x(t) = \mathop{\hbox{\rm mes}}(\{ \modo f > t \}) $.
Thus the decreasing rearrangement
$x^*(t)$ is essentially the inverse function of $d_x$.
Given functions
$x_1$, $x_2,\dots,x_n$ on $[0,1]$, we define their dilated disjoint sum
to be the function on $[0,1]$:
$$ C(x_1,\dots,x_n)(t) =
   \modo{x_k(nt-k+1)} \quad ((k-1)/n < t < k/n) .$$
Thus
$$ d_{C(x_1,\dots,x_n)}(t) = {1\over n} \sum_{k=1}^n d_{x_k}(t) .$$

We will say that an r.i.\ space $E$\ is $D$-convex if there is a constant
$c>0$ such that
$$ \normo{C(x_1,\dots,x_n)}_E \le c \, \sup_{1 \le k \le n} \normo{x_k}_E ,$$
and that $E$ is $D^*$-convex if there is a constant
$c > 0$ such that
$$ \normo{C(x_1,\dots,x_n)}_E \ge c^{-1} \,
   \inf_{1 \le k \le n} \normo{x_k}_E .$$

There is another way to define these notions.  Let us consider the
vector space $V$ of right continuous
functions from $[0,\infty)$ to $\R$ of bounded
variation.  Define the subsets
$$ \eqalignno{
   B_c^\le &= \{ d_x : \normo x_E \le c \} ,\cr
   B_c^\ge &= \{ d_x : \normo x_E \ge c \} ,\cr
   B_c^= &= \{ d_x : \normo x_E = c \} .\cr}$$
Then
$E$ is
$D$-convex if and only if
there exists a constant $c>0$ such that $\conv B_1^\le$ is
contained in $B_c^\le$, and
$E$ is
$D^*$-convex if and only if
there exists a constant $c>0$ such that $\conv B_1^\ge$ is
contained in $B_c^\ge$.

Note that Lorentz spaces as defined in Section~1 are all $D^*$-convex,
Marcinkie\-wicz spaces are all $D$-convex, and Orlicz spaces are both
$D$ and $D^*$-convex.  An easy argument shows that if $E$ is $D$-convex,
then $E'$ is $D^*$-convex, and it follows from Corollary~24 below that
if $E$ is $D^*$-convex, then $E'$ is $D$-convex.

Suppose that $M:[0,\infty) \to [0,\infty)$ is increasing.  We will say
that $M$ is $p$-convex if $M(t^{1/p})$ is convex, and
we will say that
$M$ is $q$-concave if $-M(t^{1/q})$ is convex.  By convention, we
will say that
$M$ is always $\infty$-concave.  We have the following result.

\proclaim{Lemma 20}  Suppose that $M:[0,\infty) \to [0,\infty)$
is such that there exist $1 \le p < q \le \infty$ and a constant $c>0$
such that for all $0<s<1$ that
$$ c^{-1} s^q M(t) \le M(st) \le c s^p M(t) ,$$
(where we shall suppose that the first inequality is missing if
$q = \infty$).
Then there exists an increasing, $p$-convex, $q$-concave function
$M_1$ such that there exists a constant $c_1>0$ with
$c_1^{-1} M \le M_1 \le c_1 M$.
\endproclaim

\demo{Proof}
Let $M_2(t) = \sup_{s<1} M(st)/s^p$, and let
$M_3(t) = \inf_{s<1} M(st)/s^q$\ ($M_3 \allowmathbreak
 = M_2$ if $q = \infty$).
 From now on, if $q = \infty$, we shall suppose that any
inequality involving
$q$ is automatically true.
Then
$c^{-1} M \le M_3 \le cM$, and
$$ s^q M_3(t) \le M_3(st) \le s^p M_3(t) ,$$
that is,
$M_3(t^{1/p})/t$ is an increasing function, and
$M_3(t^{1/q})/t$ is a decreasing function.  Now set
$$ \eqalignno{
   M_1(t)
   &=
   \int_0^t {M_3(s) \over s} \, ds \cr
   &=
   \int_0^{t^p} {M_3(s^{1/p}) \over p s} \, ds \cr
   &=
   \int_0^{t^q} {M_3(s^{1/q}) \over q s} \, ds ,\cr }$$
where the last equality holds only if $q<\infty$.
Then $M_1$ is $p$-convex and $q$-concave.  Further,
$M_1 \le M_3/p$, and
$$ M_1(t)
   \ge
   \int_{t^p/2}^{t^p} {M_3(s^{1/p}) \over p s} \, ds
   \ge
   1/(2p) M_3(t/2^{1/p})
   \ge 1/(4p) M_3(t) .$$
\qed
\enddemo

If $1 \le p < q \le \infty$, we say that $E$ is an interpolation space for
$(L_p,L_q)$ if there is a constant $c>0$ such that
whenever $T:L_p\cap L_q \to L_p \cap L_q$ is a linear
operator, such that $\normo T_{L_p \to L_p} \le 1$ and
$\normo T_{L_q \to L_q} \le 1$, then $\normo T_{E\to E} \le c$.

The following result is an immediate consequence of results in \cite{HM} and
Lemma~20 (see also \cite{AC}).

\proclaim{Theorem C}
Suppose that $E$ is an interpolation space for $(L_p,L_q)$.  Then
there is a constant $c>0$ such that
whenever $\normo x_{L_M} \le \normo y_{L_M}$ for all
increasing
$p$-convex and
$q$-concave functions $M:[0,\infty) \to [0,\infty)$
and if $y \in E$, then $x \in E$ and $\normo x_E \le c\, \normo y_E$.
\endproclaim

Now let us state the main results of this section.

\proclaim{Theorem 21}
Suppose that $E$ is a $D$-convex interpolation
space for $(L_p,L_q)$, where $1 \le p < q \le \infty$.
Then there exists a constant $c>0$
such that for every $x \in E$ with $\normo x_E = 1$, there exists
an increasing, $p$-convex, $q$-concave function $M:[0,\infty)
\to [0,\infty)$ such that $\int M(\modo x) \, ds \ge c^{-1}$, and
$\int M(\modo y) \, ds \le c$ whenever $\normo y \le c^{-1}$.

Thus there exists a family of increasing, $p$-convex, $q$-concave
functions $M_\alpha:[0,\infty) \to [0,\infty)$ ($\alpha\in A$)
such that $\normo \cdot_E$ is equivalent to
$\sup_{\alpha\in A} \normo\cdot_{L_{M_\alpha}}$.
\endproclaim

\proclaim{Theorem 22}
Suppose that $E$ is a $D^*$-convex interpolation
space for $(L_p,L_q)$, where $1 \le p < q \le \infty$.
Then there exists a constant $c>0$
such that for every $x \in E$ with $\normo x_E = 1$, there exists
an increasing, $p$-convex, $q$-concave function $M:[0,\infty)
\to [0,\infty)$ such that $\int M(\modo x) \, ds \le c$, and
$\int M(\modo y) \, ds \ge c^{-1}$ whenever $\normo y \ge c$.
\moreproclaim
Thus there exists a family of increasing, $p$-convex, $q$-concave
functions $M_\alpha:[0,\infty) \to [0,\infty)$ ($\alpha\in A$)
such that $\normo \cdot_E$ is equivalent to
$\inf_{\alpha\in A} \normo\cdot_{L_{M_\alpha}}$.
\endproclaim

\proclaim{Theorem 23}  Suppose that $E$ is $D$-convex
and $D^*$-convex.  Then there exists an increasing function
$M:[0,\infty)
\to [0,\infty)$ such that $E$ is equivalent to $L_M$.
\endproclaim

\proclaim{Corollary 24}
Suppose that $E$ is an interpolation
space for $(L_p,L_q)$, where $1 \le p < q \le \infty$.
If $E$ is $D$-convex, then $E$ is $p$-convex, and if $q < \infty$ then
there is a constant $c>0$ such that
given functions
$x_1$, $x_2,\dots,x_n$ on $[0,1]$
$$ \normo{C(x_1,\dots,x_n)}_E \le c \,
   \left({1\over n} \sum_{k=1}^n \normo{x_k}_E^q \right)^{1/q}.$$
If $E$ is $D^*$-convex, then $E$ is $q$-concave, and there is a
constant $c>0$ such that
given functions
$x_1$, $x_2,\dots,x_n$ on $[0,1]$
$$ \normo{C(x_1,\dots,x_n)}_E \ge c^{-1} \,
   \left({1\over n} \sum_{k=1}^n \normo{x_k}_E^p \right)^{1/p}.$$
\endproclaim

\demo{Proof}
Let us provide the proof of the stated inequality in the case that $E$ is
$D$-convex.  The other
results have almost identical proofs.

 From Theorem~21, we see that it is sufficient
to show that if $M:[0,\infty)\to[0,\infty)$ is increasing,
convex, and
$q$-concave (with $q<\infty$), then
given functions
$x_1$, $x_2,\dots,x_n$ on $[0,1]$, we have that
$$ \normo{C(x_1,\dots,x_n)}_{L_M} \le
   \left({1\over n} \sum_{k=1}^n \normo{x_k}_{L_M}^q \right)^{1/q}.$$
Let us suppose that the left hand side is bounded below by $1$.  Thus
$$ \int M(\modo{C(x_1,\dots,x_n)}) \, ds
   =
   {1\over n} \sum_{k=1}^n \int M(\modo{x_k}) \, ds
   \ge 1 .$$
Thus there exists a sequence $c_k \ge 0$ with $\sum_{k=1}^n c_k = n$
such that
$$ \int M(\modo{x_k}) \, ds \ge c_k .$$
Since $M$ is $q$-concave, it follows that
$$ \int M(\modo{x_k}/c_k^{1/q}) \, ds \ge 1 ,$$
that is, $\normo{x_k}_{L_M} \ge c_k^{1/q}$.  The result follows.
\qed
\enddemo

Let us now proceed with the proofs of the main theorems.

\proclaim{Lemma 25}  If $E \ne L_\infty$, then for each
$\epsilon>0$, there exists a
strictly increasing
function
$N:[0,\infty) \to [0,\infty)$ with $N(0) = 0$,
such that if $\int_0^1 N(x^*(t))\,dt \le 1$, then
$\normo x_E < \epsilon$.
\endproclaim

\demo{Proof}
Let $k \in E \setminus L_\infty$ such that
$k^*(1) \le 1$, $\normo k_E < \epsilon$, and $k^*$ is strictly decreasing.
Define
$$ N(t) = \cases 1/(k^*)^{-1}(t) & \text{ if }\ t \ge 1 \\
                     t/(k^*)^{-1}(1) & \text{ if }\ t \le 1 .\endcases$$
Now suppose that $\normo f_{L_N} \le 1$.  Then
$$ \int_0^1 N(f^*(t)) \, dt \le 1 ,$$
which implies that $tN(f^*(t)) \le 1$, that is $f^*(t) \le k^*(t)$.
Therefore $\normo f_E \le \normo k_E < \epsilon$.
\qed
\enddemo

For any $L>0$, we will write
$$ V_L = \{ f \in V : f(t) = 0 \hbox{ for $t > L$}\}. $$
Note that $d_x \in V_L$ if and only if $\normo x_\infty \le L$.  Notice
also that $V_L$ has a predual, $C([0,L])$, defined by the pairing
$$ \langle f,N \rangle = - \int_{[0,L]} N(s) df(s) .$$
Notice that if $\normo x_\infty \le L$, then
$$ \langle d_x, N \rangle = \int_0^1 N(\modo{x(s)}) \, ds .$$

\proclaim{Lemma 26}  Suppose that $E \ne L_\infty$, and that $L>0$.
Then the set $B_c^= \cap V_L$ is weak* compact in $V_L$.
\endproclaim

\demo{Proof}  It is clear that $B_c^=$ is a bounded set in $V$, and
hence it is sufficient to show that $B_c^\le \cap V_L$ is weak* closed
in $V_L$.
Suppose that $\normo{x_n}_E = c$, that $d_{x_n} \in V_L$, and that
$d_{x_n} \to g$ weak*.
Then it is easy to see that $g$\ is decreasing with $g(0) \le 1$, and that
$d_{x_n} \to g$ pointwise except possibly at discontinuities of $g$.
Therefore $g = d_y$ for
some $y \in V_L$, and
$x_n^* \to y^*$ pointwise except possibly at points of
discontinuity of $y^*$, of which there are only countably many.
Hence by Lebesgue's Dominated Convergence Theorem, it
follows that for any continuous function $N$ that
$\int_0^1 N(x_n^* - y^*) \, dt \to 0$,
and hence by Lemma~1, it follows that $\normo{x^*_n-y^*}_E
\to 0$.
\qed
\enddemo

Now, if $1 \le p < q \le \infty$, and $L>0$, we define the subset
$C_{p,q,L}$ of
$V_L$ to be the set of all those $f \in V_L$ such that
$$ -\int_{[0,L]} N \, df \ge 0 $$
for all increasing $N:[0,\infty) \to [0,\infty)$
that are $p$-convex and $q$-concave.
Notice that $C_{p,q,L}$ is weak* closed in
$V_L$ for all $L > 0$.

\bigskip

\demo{Proof of Theorem 21}
We may suppose that $E \ne L_\infty$.
Since we have that $\normo{\min\{\modo x,t\}}_E 
 \to \normo x_E$ as $t\to\infty$, we may
suppose without loss of generality that $x \in L_\infty$.  By a further
slight approximation, we may suppose that $x^*$ is strictly decreasing
and $x^*(1) = 0$, that is, we may suppose that $d_x$ is absolutely
continuous.

By the definition of $D$-convexity, and Theorem~C,
we know that there exists
a constant $c>0$ such that for all $L>0$ we have that
$\conv(B_{c^{-1}}^=)$ does not intersect with $\{d_{x}\} + C_{p,q,L}$.
Hence, by the Hahn-Banach Theorem, for each $L>\normo x_\infty$,
there exists
$M_L \in C_0([0,L])$\
such that for some constant $S = \pm 1$
$$ \int_0^1 M_L(\modo y) \, ds =
   -\int_0^L M_L \, d(d_y) \le S \qquad \hbox{for} \quad
   d_y \in B_{c^{-1}}^\le \cap V_L, $$
and
$$ -\int_0^L M_L \, d(d_{cx} + f) \ge S \qquad \hbox{for} \quad
f \in C_{p,q,L}.$$
Hence
$$ \int_0^1 M_L(\modo x) \, ds =
   -\int_0^L M_L d(d_x) \ge S .$$
Further, since $C_{p,q,L}$ is a cone, it follows that
$$ -\int_0^L M_L \, df \ge 0 \qquad \hbox{for} \quad
f\in C_{p,q,L}.$$
For $0\le a < b \le L$, consider the functions
$$ f_1(s) = \cases
   1  & \text{ if }\ a \le s \le b \\
   0  & \text{ otherwise}, \endcases $$
$$ f_2(s) = \cases
   -1 & \text{ if }\ a^{1/p} \le s \le ((a+b)/2)^{1/p} \\
   1  & \text{ if }\ ((a+b)/2)^{1/p} \le s \le b^{1/p} \\
   0  & \text{ otherwise},\endcases $$
$$ f_3(s) = \cases
   1  & \text{ if }\ a^{1/q} \le s \le ((a+b)/2)^{1/q} \\
   -1 & \text{ if }\ ((a+b)/2)^{1/q} \le s \le b^{1/q} \\
   0  & \text{ otherwise}.\endcases $$
It is easily seen that $f_1$, $f_2$ and $f_3$ are in $C_{p,q,L}$,
and hence it may be seen that
$M_L$ is positive, increasing, $p$-convex, and $q$-concave on $[0,L]$.
Hence $S=1$.

By Lemma~25, there exists $N:[0,\infty) \to [0,\infty)$
with $N(0) = 0$
and that is strictly increasing, and such that if $\int_0^1
N(\modo{z}) \, ds \le 1$,
then $\normo z_E \le 1$.  Hence, if $z$ with $\normo z_\infty \le L$, and
if $\int_0^1 N(\modo{z}) \, ds \le 1$, then $\int M_L(\modo{z}) \, ds
\le 1$.  Hence $M_L \le N$.

Notice that $(M_L/N)$ is a bounded sequence
in $L_\infty([0,\infty))$.  Let $M/N$ be a weak* limit point of this
sequence.  Since $N(s) {d\over ds}(d_z(s))$ is in
$L_1([0,\infty))$ whenever $z \in L_N$, and
$d_x$ is absolutely continuous,
it is easy to see
that $M$ satisfies the requirements of Theorem~1.  (Initially one
would have to restrict to those $y \in L_\infty$, but an application
of Lebesgue's Monotone Convergence Theorem will deal with this.)
\qed
\enddemo

\demo{Proof of Theorem~22}
Following the first part of the proof of Theorem~21 with only minor
modifications, we can show the following.  If $x \in L_\infty$
with $\normo x_\infty = L$, then there exists an increasing,
$p$-convex, $q$-concave function $M_1:[0,L]\to[0,\infty)$
such that
$$ \int M_1(\modo x) \, ds \le c ,$$
and whenever $\normo y_\infty \le L$ with $\normo y_E \ge c$, then
$$ \int M_1(\modo y) \, ds \ge 1 .$$

If $q = \infty$, consider the function $N$ generated by Lemma~25 in
the case when $\epsilon = 1$.
Notice that if we set
$$ N_1(t) = \cases
   N(t) & \text{ if }\ t \ge 1 \\
   N(1) t^p & \text{ if }\ t<1, \endcases$$
and
$$ N_2(t) = \sup_{s<1} N_1(st)/s^p ,$$
then $N_2$ still satisfies the conclusion of Lemma~25.  Furthermore, if
we set
$$ M(t) = \cases
   M_1(t) & \text{ if }\ t \le L \\
   (M_1(L)/N_2(L))N_2(t) & \text{ if }\ t > L,\endcases$$
then we see that $M$ satisfies the hypotheses of Lemma~20.  Also,
if $\normo y_E \ge 2c$, with $y \ge 0$, then either
$\normo{\min\{y,L\}}_E
\ge c$, in which case
$\int M(\min\{y,L\}) \, ds = \int M_1(\min\{
y,L\}) \, ds \ge c^{-1}$,
or $\normo{yI_{y>L}}_E \ge c$, in which case
$\int M(yI_{y>L}) \, ds = (M_1(L)/N_2(L)) \int N_2
(yI_{y>L}) \, ds \ge c^{-1}$.  In either case,
$\int M(y) \, ds \ge c^{-1} $.

If $q < \infty$, it is an easy matter to see
that $L_q \subseteq E$, and that there exists a constant $c_1>0$
such that $\normo z_E \le c_1 \normo z_q$ for all $z \in L_q$.
Set
$$ M(t) = \cases
   M_1(t) & \text{ if }\ t \le L \\
   M_1(L) (t/L)^q & \text{ if }\ t > L.\endcases$$
Thus $M$ satisfies the hypotheses of Lemma~20.  Furthermore, if
$\normo y_E \ge 2c$, then
$\int M(\modo y) \, ds \ge c^{-1}$, by the same argument as in the
case when $q = \infty$.

Now let us consider the case for general $x$.  Note
that $E \subseteq L_p$, and that there is a constant
$c_2>0$ such that $\normo z_p \le c_2 \normo z_E$ for all $z \in E$.
Without loss of generality, $x\ge 0$.
Write $x = x_1 + x_2$, where
$x_1$ and $x_2$ have disjoint support,
$x_1 \in L_\infty$, $\normo{x_1}_E \ge 1/2$,
and $\normo{x_2}_p \le 1$.

Let $M_1$ be the function described by Theorem~22 for $x_1$, and let
$M(t) = \min\{M_1(t),t^p\}$.  It is clear that $M$ satisfies
the hypotheses of Lemma~20, and also that $\int M(x) \, ds \le c$.
Now suppose that $\normo y \ge 4c$.  We may suppose that $y \ge 0$.
Write $y = y_1 + y_2$, where $y_1 = y I_{M_1(y) \le y^p}$.  Then
either $\normo {y_1}_E \ge 2c$, in which case
$\int M(y) \, ds \ge \int y_1^p \, ds \ge 2c_2^{-1}$, or
$\normo{y_2}_E \ge 2c$, in which case
$\int M(y) \, ds \ge \int M_1(y_2) \, ds \ge c^{-1}$.
\qed
\enddemo

We will leave the proof of Theorem~23 to the reader, as it follows the
ideas of the previous proofs.

We also leave with a problem that was given to the first named author
by
Carsten Sch\"utt.
If $E$ is $D^*$-convex, does there exist an
appropriately measurable
family
of increasing, convex functions
$M_\alpha:[0,\infty) \to [0,\infty)$ ($\alpha\in A$), where $A$
is a measurable space with measure $\mu$,
such that
$\normo\cdot_E$ is equivalent to
$\int \normo\cdot_{M_\alpha} \, d\mu(\alpha)$?

\head{7.  Another generalization of Theorem~B}\endhead
In this section we will consider another generalization of Theorem~B.
Suppose that $X$ is a symmetric sequence space on sequences
$x = (x_i)_{1\le i \le n} \in \R^n$.
Let us suppose that $\normo{(1,0,\dots,0)}_X = 1$.
Then we define its associated
r.i.\ space, $E_X$ by the following formula:
$$ \normo x_{E_X}
   =
   \normo{ \left( {1\over n} \int_{(i-1)/n}^{i/n} x^*(s) \, ds \right)
   _{1\le i \le n} }_X .$$
Let us show that $E_X$ really does satisfy the triangle inequality.  It
is clear that $\normo{x^*+y^*}_{E_X} \le \normo x_{E_X} + \normo y_{E_X}$.
It is also easy to see that if $x\prec y$, then $\normo x_{E_X} \le
\normo y_{E_X}$.  Since $x+y \prec x^* + y^*$, we are done.

To save space, if $A$ and $B$ are two quantities
depending upon certain parameters, we will write
$A \approx B$ if there exists a constant $c>0$, independent of the
parameters, such that $c^{-1} A \le B \le c A$.  If $t$ is a real number,
we will write $[t]$ for the greatest integer less than $t$.

\proclaim{Theorem~27}
There exists a constant $c>0$
such that
if $(x_{i,j})_{1\le i,j \le n}$ is an $n\times
n$ matrix, then
$$ {1\over n!} \sum_{\pi \in S_n}
   \normo{(x_{i,\pi(i)})_{1\le i \le n}}_X
   \ge
   c^{-1} \, \left(
   {1\over n} \sum_{k=1}^n x^*_k
   + \normo{(x^*_{kn})_{1\le k\le n})}_X
   \right) .$$
Furthermore, if the associated r.i.\ space is $D^*$-convex, then
there exists a constant $c>0$
such that
if $(x_{i,j})_{1\le i,j \le n}$ is an $n\times
n$ matrix, then
$$ {1\over n!} \sum_{\pi \in S_n}
   \normo{(x_{i,\pi(i)})_{1\le i \le n}}_X
   \le
   c \, \left(
   {1\over n} \sum_{k=1}^n x^*_k
   + \normo{(x^*_{kn})_{1\le k\le n})}_X
   \right) .$$
\endproclaim

We do not know whether the condition that the associated space be
$D^*$-convex is necessary in order for the second inequality to hold.
In order to show this result, we will use the following result due to
Kwapie\'n and Sch\"utt.

\proclaim{Theorem~D {\rm \cite{KS2}}}  There exist a constant $c>0$ such that
for any $n\times n\times n$ array $y = (y_{i,j,k})_{1\le i,j,k \le n}$,
we have that
$$ {1\over (n!)^2} \sum_{\pi,\sigma \in S_n}
   \max_{1 \le i \le n} \modo{ y_{i,\pi(i),\sigma(i)} }
   \approx
   {1 \over n^2} \sum_{k=1}^{n^2} y_k^* .$$
\endproclaim

\demo{Proof of Theorem~27}  Let us first consider the case when
$X_m$ is the symmetric sequence space given by
$$ \normo z_{X_m} = \sum_{k=1}^m z^*_k .$$
Suppose that given $z$, one forms the array
$$ y_{i,j} = \cases
   z_i & \text{ if }\ j \le n/m \\
   0 & \text{ otherwise}. \endcases $$
Then by Theorem~A, it may be seen that
$$ \eqalignno{
   {1\over n!} \sum_{\pi \in S_n} \max_{1\le i \le n}
   \modo{y_{i,\pi(i)}}
   &\approx
   {1\over n} \sum_{k=1}^n y_k^* \cr
   &\approx
   {[n/m] \over n} \sum_{k=1}^{n/[n/m]} z_k^* \cr
   &\approx
   {1\over m} \normo z_{X_m} . \cr } $$
Now, given $x$ as in the hypothesis of the theorem,
form the following array:
$$ y_{i,j,k} = \cases
   x_{i,j} & \text{ if }\ k \le n/m \\
   0 & \text{ otherwise}. \endcases $$
In that case
$$ \eqalignno{
   {1\over n!} \sum_{\pi \in S_n}
   \normo{(x_{i,\pi(i)})_{1\le i \le n}}_{X_m}
   &\approx
   {m\over (n!)^2} \sum_{\pi,\sigma \in S_n}
   \max_{1 \le i \le n} \modo{ y_{i,\pi(i),\sigma(i)} } \cr
   &\approx
   {m \over n^2} \sum_{k=1}^{n^2} y_k^* \cr
   &\approx
   {m[n/m] \over n^2} \sum_{k=1}^{n^2/[n/m]} x_k^* \cr
   &\approx
   {1\over n} \sum_{k=1}^{nm} x_k^* \cr
   &\approx
   {1\over n} \sum_{k=1}^n x_k^*
   +
   \normo{(x_{kn}^*)_{1 \le k \le n}}_{X_m} .\cr }$$
Now let us consider more general symmetric sequence spaces $X$.
We know that
$$ \eqalignno{
   \normo z_X
   &=
   \sup_{\normo w_{X^*} \le 1}
   \sum_{k=1}^n z_k^* w_k^*  \cr
   &=
   \sup_{\normo w_{X^*} \le 1}
   \sum_{m=1}^n (w_m^*-w_{m+1}^*) \normo z_{X_m} ,\cr }$$
where by convention $w^*_{n+1} = 0$.  From this, we immediately see
that for some constant $c>0$
$$ \eqalignno{
   {1\over n!} \sum_{\pi \in S_n}
   \normo{(x_{i,\pi(i)})_{1\le i \le n}}_X
   &\ge c^{-1}\,
   \sup_{\normo w_{X^*} \le 1}
   \sum_{m=1}^n (w_m^*-w_{m+1}^*)
   \left(
   {1\over n} \sum_{k=1}^n x_k^*
   +
   \normo{(x_{kn}^*)_{1 \le k \le n}}_{X_m}
   \right)  \cr
   &\ge c^{-1}\,
   \left(
   {1\over n} \sum_{k=1}^n x_k^*
   +
   \normo{(x_{kn}^*)_{1 \le k \le n}}_X
   \right) . \cr}$$
since $w_1^* \le 1$ whenever $\normo w_{X^*} \le 1$.

Now let us show the second inequality when $E_X$ is $D^*$-convex.
Let us consider the following functions:
$$ \eqalignno{
   z_\pi(t) &=
   x_{i,\pi(i)} \quad t \in [(i-1)/n,i/n) , \cr
   w(t) &= x_k^* \quad t \in [(k-1)/n^2,k/n^2) . \cr } $$
It is an easy matter to see that
$ w^* = C(z_\pi:\pi\in S_n)^* $.
Hence, by Corollary~24, we see that for some constant $c>0$ depending
only on $X$
$$ {1\over n!} \sum_{\pi \in S_n} \normo{z_\pi}_{E_X}
   \le
   c \, \normo w_{E_X} .$$
The result now follows after we notice that
$$ \normo x_{E_X}
   \approx
   {1\over n} \int_0^{1/n} x^*(s) \, ds
   +
   \normo{(x^*(k/n))_{1\le k \le n}}_X .$$
\qed
\enddemo

\head{8. $D$ and $D^*$-convex Lorentz Spaces}\endhead
Although the results in this section are primarily concerned with
Lorentz spaces, in order to prove our results, we will need a wider
class of spaces, known as Orlicz-Lorentz
spaces.  If $M,N:[0,\infty) \to [0,\infty)$ are strictly increasing
bijections, then we define the space $L_{M,N}$ to be the set of those
measurable
functions $x:[0,1]\to \R$ such that
$$ \normo x_{L_{M,N}}
   = \normo{x^*\circ \tilde M \circ \tilde N^{-1}}_{L_N} ,$$
where $\tilde M = 1/M(1/t)$, and $\circ$ denotes function composition.
It is not clear what are necessary and sufficient conditions for
$L_{M,N}$ to have an equivalent norm that satisfies the triangle inequality,
but this will not be relevant to our discussion.  It is clear that
$L_M = L_{M,M}$, and that $\Lambda_r(\varphi) = L_{\tilde \varphi^{-1},t^r}$
with equality of norms.

Following \cite{M1}, we say that an increasing bijection $M:[0,\infty) \to
[0,\infty)$ is almost convex if there are numbers
$a>1$, $b>1$, and a positive integer $p$ such that for all positive integers
$m$, the cardinality of the set of integers $n$ such that we do not have
$M(a^{n+m}) \ge a^{m-p} M(a^n)$
is less than $b^m$.  It is clear that this notion also can be made to
make sense if
$M$ is only a bijection from $[0,1] \to [0,1]$, or a bijection from
$[1,\infty) \to [1,\infty)$, by stating that the inequality is true
whenever it is undefined.

The following result is essentially Theorem~4.2 from \cite{M1}.  The results
from \cite{M1} are concerned with function spaces on $\R$ rather than
$[0,1]$, but the change is not too hard to do.

\proclaim{Theorem E}  Let $M$, $N_1$, $N_2$ be increasing
bijections $[0,\infty)\to[0,\infty)$ that map $1$ to $1$,
such that one of $N_1$ or $N_2$ is convex and $q$-concave for some
$q<\infty$.
Then the following are
equivalent.
\roster\widestnumber\item{ii)}
\item"{i)}" For some $c>0$ we have that
$\normo x_{L_{M,N_1}} \le c \, \normo x_{L_{M,N_2}}$ for all
measurable $x:[0,1] \to \R$.
\item"{ii)}" $N_1 \circ N_2^{-1}$ restricted to $[1,\infty)$ is almost
convex.
\endroster
\endproclaim

\proclaim{Theorem 28}  If $\varphi:[0,1] \to [0,1]$ is an increasing
bijection, and $1 \le r < \infty$, such that $\Lambda_r(\varphi)$ is
equivalent to a norm, then $\Lambda_r(\varphi)$ is $D$-convex if and
only if $(\tilde \varphi(t))^r$ is almost convex, and $D^*$-convex
if and only if
$\tilde\varphi^{-1}(t^{1/r})$ is almost convex.
\endproclaim

\demo{Proof}
Suppose that $\Lambda_r(\varphi)$ is $D$-convex.  Define
$$ M(t) = \cases t & \text{ if }\ 0\le t < 1 \\
                  \tilde \varphi^{-1}(t) & \text{ if }\ t \ge 1 . 
\endcases $$
It is clear that $\normo{I_{[0,t]}}_{\Lambda_r(\varphi)} =
\normo{I_{[0,t]}}_{L_M} = \varphi(t)$ for $0\le t \le 1$, and that
$ \Lambda_r(\varphi) = L_{M,t^r} $.

We will show that there is a constant $c_1>0$ such that
$\normo x_{\Lambda_r(\varphi)} \le c_1 \normo x_{L_M}$.  Then the
result will follow easily from Theorem~E.

For, by Theorem~21, there exists a constant $c_2>0$ so that the following
holds.
Suppose that $\normo x_{\Lambda_r(\varphi)} = 1$.  Then
there exists an increasing convex bijection $N:[0,\infty) \to [0,\infty)$
such that $\normo x_{L_N} \ge c_2^{-1}$, but that in
general
$\normo y_{L_N} \le \normo y_{\Lambda_r(\varphi)}$.  By considering
$y = I_{[0,t]}$, we see that $N(t) \le M(t)$ for $t \ge 1$.  Now,
for any $\lambda < c_2^{-1}$, we have that
$$ \int_0^1 N(\modo x/\lambda) \, ds > 1 .$$
Since $N$ is convex, for any $\lambda < c_2^{-1}/2$, we have that
$$ \int_0^1 N(\modo x/\lambda) \, ds > 2.$$
Further, for any $\lambda > 0$
$$ \int_{\modo x \le \lambda} N(\modo x/\lambda) \, ds \le N(1) \le 1 .$$
Hence for $\lambda < c_2^{-1}/2$, we have that
$$ \int_0^1 M(\modo x/\lambda) \, ds \ge
   \int_{\modo x \ge \lambda} M(\modo x/\lambda) \, ds \ge
   \int_{\modo x \ge \lambda} N(\modo x/\lambda) \, ds > 1 , $$
and hence $\normo x_{L_M} \ge c_2^{-1}/2$.

The case when $\Lambda_r(\varphi)$ is $D^*$-convex is almost identical.
\qed
\enddemo

\enddocument